\magnification=\magstep1
\input amstex
\documentstyle{amsppt}
\voffset=-3pc
\loadbold
\loadmsbm
\loadeufm
\UseAMSsymbols
\baselineskip=12pt
\parskip=6pt
\def\var{\varepsilon}

\def\bP{\Bbb P}

\def\uU{\frak U}
\def\uV{\frak V}
\def\uW{\frak W}

\def\bC{\Bbb C}
\def\bN{\Bbb N}
\def\Ker{\text{ Ker}}
\def\id{\text{ id}}
\def\Hom{\text{Hom}}
\def\cF{\Cal F}
\def\cO{\Cal O}
\NoBlackBoxes
\topmatter
\title On the Cohomology Groups of Holomorphic Banach Bundles\endtitle
\rightheadtext{Cohomology Groups of Banach Bundles}
\address Department of Mathematics,
Purdue University,
West Lafayette, IN\ 47907-1395\endaddress
\author
L\'aszl\'o Lempert\endauthor\footnote""{Research partially supported by NSF
grants DMS0203072 and 0700281\hfill\break}
\abstract
We consider a compact complex manifold $M$, and introduce the notion of two holomorphic Banach bundles
$E,F$ over $M$ being compact perturbations of one another.
Given two such bundles we show that if the cohomology groups $H^q(M,E)$ are finite dimensional then so
are the cohomology groups $H^q(M,F)$; as well as a more precise result in the same spirit.
\endabstract
\endtopmatter
\document
\head Introduction\endhead

A basic result in complex geometry is the finiteness theorem of Cartan and Serre [CS], according to which the
cohomology groups of a finite rank holomorphic vector bundle $E\to M$ over a compact base are finite
dimensional.
(More generally, the same holds for the cohomology groups of coherent sheaves over $M$.)
It is no surprise that holomorphic {\sl Banach} bundles $E\to M$ over a compact base can easily have infinite
dimensional cohomology groups.
For example, J.~Kim observes that if $E$ is the trivial bundle with fiber a Banach space $A$, then
$H^q(M,E)\approx H^q(M,\Cal O)\otimes A$ (tensor product over $\bC$); the special case when $H^q(M,\Cal O)=0$
is due to Leiterer [K, Li2].
In particular, $\dim H^q (M,E)=\infty$ if $\dim A=\infty$ and $H^q(M,\Cal O)\neq 0$.

However, already in [G] Gohberg was led to a class of holomorphic (first Hilbert, then) Banach bundles, for
which finiteness theorems can be proved.
The terminology is {\sl bundles of compact type}, meaning their structure group can be reduced to the group of
invertible operators of form id + compact.
Let $F\to\bP_1$ be a holomorphic Banach bundle of compact type over the Riemann sphere.
The Riemann--Hilbert type result of [G,GL], coupled with Bungart's work [B] implies that $F\approx T\oplus V$, where $T$ is trivial and $V$ has finite rank.
Hence
$$
H^1(\bP_1,F)\approx H^1 (\bP_1,T)\oplus H^1 (\bP_1,V)
\approx H^1(\bP_1,V)
$$
is finite dimensional by what has been said above.
More recently Leiterer in [Li2] considers an arbitrary compact complex manifold $M$, and holomorphic Banach
bundles $F\to M$ of compact type, $V\to M$ of finite rank.
Assuming $H^q(M,V)=0$ for some $q\geq 0$, and putting a mild restriction on the fiber type of $F$, he proves $\dim H^q (M,F\otimes V) < \infty$.

Our goal in this paper is to generalize and thereby better to understand the phenomenon discovered by Gohberg
and Leiterer.
Instead of bundles of compact type, we shall work with a bundle and its compact perturbation, by which we mean
the following.
Let $N$ be a complex manifold and $E,F\to N$ holomorphic Banach bundles.

\definition{Definition 1.1}(a) A fredomorphism $E\to F$ is a homomorphism $\varphi\colon E\to F$ such that $\varphi|E_x:E_x\to
F_x$ is a Fredholm operator for each $x\in N$.

(b) We say that $F$ is a compact perturbation of $E$ if there are an open cover $\frak U$ of $N$ and fredomorphisms
$\varphi_U\colon E|U\to F|U$ such that 
$\varphi_U-\varphi_V|E_x\colon E_x\to F_x$ is compact whenever 
$U,V\in\frak U$ and $x\in U\cap V$.
\enddefinition

It is easy to see (cf.~Proposition 4.2) that being a compact perturbation is a symmetric and transitive relation.---For example, any two finite rank bundles $E,F\to N$ are compact perturbations of each other, the zero homomorphism $E\to F$ being a (global) fredomorphism.
Further, if $F$ is a holomorphic Banach bundle of compact type whose
fibers are isomorphic to a Banach space $A$, then $F$ is a compact perturbation of the trivial bundle $A\times N\to N$.

We recall that the sheaf cohomology groups $H^q(N,E)$ of holomorphic Banach bundles are topological vector spaces, see Section 2.

\proclaim{Theorem~1.2}Let $M$ be a compact complex manifold, $E,F$ holomorphic Banach bundles over $M$ that are compact perturbations of one another, and $q=0,1,\ldots$.
If $H^{q+1}(M,E)$ is Hausdorff (in particular if $\dim H^{q+1}(M,E)<\infty$) and $\dim H^q (M,E)<\infty$, then $H^q(M,F)$ is also finite dimensional (and Hausdorff).
\endproclaim

The assumption that $H^{q+1}(M,E)$ be Hausdorff is not automatically satisfied, and without it the conclusion would not neccessarily hold,
see the discussion at the end of the paper.

There should be a sheaf version of this theorem, involving cohesive sheaves from [LP] or Leiterer's Banach coherent analytic Fr\'echet sheaves from [Li1], but so far we have only been able to prove an extension of Theorem 1.2 to certain very special subsheaves of $E,F$ above.

Compact perturbations can be defined more generally for topological
Banach bundles. If $M$ is a compact topological space and $E,F\to M$ are
Banach bundles, any system of fredomorphisms
$\varphi_U:E|U\to F|U$, $U\in\frak U$ as in Definition 1.1b induces
a virtual vector bundle of finite rank over $M$, i.e., an element
of the group $K(M)$. In particular, one can speak of characteristic classes
in this situation. Perhaps even in the holomorphic category one should
regard compact perturbations as generalizations of holomorphic vector
bundles of finite rank; Theorem 1.2 seems to point in this direction.

The finiteness results of Gohberg and Leiterer are special cases of Theorem 1.2.
Indeed, as said, a bundle $F$ of compact type is a compact perturbation of a trivial bundle $E=A\times M\to M$, and $F\otimes V$ is a compact perturbation of $E\otimes V$, when rk $V<\infty$.
Now in [K] J.~Kim proves $H^p(M,E\otimes V)\approx H^p(M,V)\otimes A$ as topological vector spaces; in particular all $H^p(M,E\otimes V)$ are Hausdorff.
(In general there are many natural ways to topologize the tensor product $Y\otimes A$ of Banach spaces; but there is no ambiguity when $\dim Y=n <\infty$ as above.
In this case $Y\otimes A\approx A\oplus\ldots\oplus A$, $n$ times.)
The Cartan--Serre theorem, at least for bundles over manifolds, is also a special case.

In analysis and geometry finite dimensionality results typically depend on functional analysis:\ Fredholm theory or Schwartz's theorem on compact perturbations of Fr\'echet space operators.
We shall also help ourselves to these tools, but the present proof will have a nonlinear ingredient as well, Michael's selection theorem [M], whose
application in complex geometry is pretty much the only novelty in this 
work.

\head 2.\ Holomorphic Fr\'echet bundles and their cohomology groups. Generalities\endhead

By a Fr\'echet space we shall mean here a complete, metrizable locally convex topological vector space; a sequence of seminorms defining the topology need not be specified.
Similarly, by a Banach space we mean a completely normable locally convex space.
A continuous linear map between locally convex spaces will be called a homomorphism.
A (locally trivial) holomorphic Fr\'echet bundle over a finite dimensional complex manifold $M$ is a holomorphic map $\pi\colon E\to M$ of a complex 
Fr\'echet manifold $E$ on $M$, with each fiber $E_x=\pi^{-1}(x)$, $x\in M$, endowed with the structure of a complex vector space.
It is required that for every $x\in M$ there be a neighborhood $U\subset M$ of $x$, a complex Fr\'echet space $A$, and a biholomorphism (trivialization) $E|U=\pi^{-1} U\to A\times U$ mapping the fibers $E_y$ linearly on $A\times \{y\}\cong A$, $y\in U$.
If $A$ above is a Banach space, we speak of a Banach bundle.
For more about these notions, see [Lm1, Sections 1,2].
We fix a holomorphic Fr\'echet bundle $E\to M$ and in this Section we describe basic properties of its sheaf cohomology groups $H^q(M,E)$ as locally convex spaces. In addition to notation and definitions introduced
in this
Section, the results that we will need later on are  Lemma 2.5
and Corollary 2.6.

Let $\varphi_U\colon E|U\to A\times U$ be a trivialization.
For each compact $K\subset U$ and continuous seminorm $p$ on $A$ we define seminorms
$$
\|e\|_{K,p}=\sup_{x\in K} p(\varphi_U e(x)),\qquad e\in\Gamma(U,E),
$$
on the space $\Gamma(U,E)$ of holomorphic sections.
They endow $\Gamma(U,E)$ with the structure of a complete locally convex space, in fact a Fr\'echet space, since a countable collection of $K,p$ suffices to define the topology.
Clearly, this Fr\'echet space structure is independent of the choice of the trivialization.

If $q=0,1,\ldots$, and $\uU$ is a countable open cover of $M$ such that $E|U$ is trivial for each $U\in \uU$, then the space of (not necessarily alternating) $q$--cochains
$$
C^q(\uU,E)=\prod_{U_0,\ldots, U_q\in\uU}
\Gamma\bigl(\bigcap^q_{i=0} U_i,E\bigr)
$$
is again a Fr\'echet space; we use the convention that the
space of functions (and later, forms) on the empty set
is the $0$ dimensional space $(0)$. The coboundary operator
$$
\delta=\delta^q_{\uU}: C^q(\uU,E)\to C^{q+1} (\uU,E)
$$
is continuous, $Z^q(\uU,E)=\Ker\delta^q_{\uU}$ is a closed subspace, and the cohomology groups
$$
H^q(\uU,E)=Z^q (\uU,E)/\delta^{q-1}_{\uU} C^{q-1} (\uU,E)
$$
are locally convex vector spaces, complete and semimetrizable, but possibly non Hausdorff.
A refinement map $\uV\to\uU$ between two covers induces homomorphisms
$\rho_\frak V^\frak U\colon H^q (\uU,E)\to H^q (\uV,E)$.
We endow the vector space direct limit $H^q(M,E)$ of the system 
$\{H^q (\uU,E),\rho_\uV^\uU\}$ with the finest locally convex topology for which the canonical maps $H^q(\uU,E)\to H^q(M,E)$ are continuous. (It can be seen from Theorem 2.2 below that this
is simply the finest topology which makes the canonical maps continuous.)

A variant of this construction provides better control of the topology obtained.
If $\varphi_U\colon E|U\to A\times U$ is a trivialization as above and 
$W\subset\subset U$ is open, then for each continuous seminorm $p$ on $A$ 
and $e\in\Gamma (W,E)$ we let 
$\|e\|_p=\sup_{x\in W} p(\varphi_U e(x))\leq\infty$, and define 
the space of bounded holomorphic sections by 
$$
\Gamma_b (W,E)=\{e\in\Gamma (W,E)\colon \|e\|_p<\infty\text{ for each } p\}.
$$
This space, with the seminorms $\|\ \ \|_p$ is a Fr\'echet space, whose topology is again independent of the choice of $U\supset\supset W$ and the trivialization $\varphi_U$.
Obviously, the embedding $\Gamma_b (W,E)\to\Gamma(W,E)$ is continuous.
If $\uW$ is a countable open cover of $M$ consisting of such $W$, then the space of bounded $q$--cochains,
$$
C_b^q (\uW,E)=\prod_{W_0,\ldots,W_q\in\uW}
\Gamma_b\bigl(\bigcap^q_{i=0} W_i,E\bigr)
$$
is a Fr\'echet space, and the coboundary $\delta$ restricts to a homomorphism $C_b^q(\uW,E)\to C_b^{q+1}(\uW,E)$.
Denoting the kernel of $\delta|C^q_b(\uW,E)$ by
$Z_b^q(\uW,E)\subset C_b^q(\uW,E)$, the bounded cohomology groups
$$
H_b^q (\uW,E)=Z_b^q (\uW,E)/\delta_{\uW}^{q-1} C^{q-1} (\uW,E)
$$
are complete semimetrizable locally convex spaces.
With a cover $\uU$ finer than $\uW$, restriction induces homomorphisms
$$
H_b^q (\uW,E)\to H^q (\uU,E).\tag2.1
$$
Composing these with the canonical homomorphisms $H^q(\uU,E)\to H^q (M,E)$ we obtain canonical homomorphisms
$$
H_b^q (\uW,E)\to H^q (M,E).\tag2.2
$$

\definition{Definition 2.1}We say that a countable open cover $\uW$ of $M$ is special (with respect to $E$) if (a) for every $W\in\uW$ there is a biholomorphism of a neighborhood of $\overline W$ into some $\bC^n$ that maps $W$ on a bounded, strongly pseudoconvex domain with smooth boundary; also $E$ is trivial on this neighborhood, and (b) the boundaries of $W\in\uW$ are in general position.
This latter means that if $k\in\bN$ and $\rho_i$ are smooth defining functions of $W_i\in\uW$, $i=1,\ldots,k$, then $d\rho_1,\ldots,d\rho_k$ are linearly independent at each point of the set $\{\rho_1=\ldots=\rho_k=0\}$.
\enddefinition

Sard's lemma implies that $M$ has arbitrarily fine special covers.

\proclaim{Theorem 2.2}If $\uW$ is a special cover of $M$, then (2.2) is an isomorphism; if, in addition, $\uU$ consists of Stein open sets then (2.1) is also an isomorphism.
\endproclaim

The proof uses routine sheaf theory plus certain deeper facts of several complex variables.
We start by introducing a few more spaces generalizing $\Gamma(U,E)$, $\Gamma_b(W,E)$, $C^q(\uU,E)$ etc.
Fix smooth vector fields $\xi_1,\ldots,\xi_m$ that span each tangent space $T_x M$.
Let $\varphi_U\colon E|U\to A\times U$ be a trivialization over an open $U\subset M$.
We endow the space $\Omega^{0,r}(U,E)$ of smooth $E$--valued $(0,r)$ forms with seminorms $\|\ \ \|_{K,p,j}$, one for each compact $K\subset U$, continuous seminorm $p$ on $A$, and $j=0,1,2,\ldots$:
$$
\|f\|_{K,p,j}=
\max\max_K p\left( \xi_{i_1}\xi_{i_2}\ldots\xi_{i_j}
(\varphi_U f(\xi_{i_{j+1}},\ldots,\xi_{i_{j+r}}))\right),\tag2.3
$$
the outer $\max$ taken over all tuples 
$1\leq i_1,\ldots,i_{j+r}\leq m$.
If $\uU$ is a countable cover of $M$ consisting of such open sets, we define a double complex
$$
C^{qr}(\uU)=\prod_{U_0,\ldots,U_q\in\uU}
\Omega^{0,r}\bigl(\bigcap^q_{i=0} U_i,E\bigr),\quad q,r\geq 0,
$$
of Fr\'echet spaces.
We let
$$
\delta=\delta^{qr}\colon C^{qr} (\uU)\to C^{q+1,r}(\uU),\qquad\overline\partial=\overline\partial^{qr}\colon C^{qr}(\uU)\to C^{q,r+1}(\uU)
$$
denote \v Cech coboundary, resp.~the homomorphism obtained by applying $\overline\partial$ componentwise.
Set
$$
Z^{qr}(\uU)=\Ker \overline\partial^{qr}\cap\Ker\delta^{qr}.
$$

With $U$ as above and $W\subset\subset U$ open, we also consider the space $\Omega_b^{0,r}(W,E)\subset\Omega^{0,r}(W,E)$ of bounded $(0,r)$ forms $f$, meaning that for each continuous seminorm $p$ on $A$ $\|f\|_p=\max\sup_W p(\varphi_U f(\xi_{i_1},\ldots,\xi_{i_r})) 
<\infty$, where the $\max$ is taken over all 
$1\le i_1,\ldots,i_r\le m$.
We endow $\Omega^{0,r}_b (W,E)$ with these seminorms $\|\ \ \|_p$ to obtain 
a metrizable locally convex space, no longer complete.
Observe that the embedding $\Omega_b^{0,r}(W,E)\to\Omega^{0,r}(W,E)$ is not continuous when $r>0$.
Nevertheless, given a countable cover $\uW$ consisting of such $W$ there are the spaces
$$
C_b^{qr}(\uW)=\prod_{W_0,\ldots,W_q\in\uW} 
\Omega_b^{0,r}\bigl(\bigcap^q_{i=0} W_i,E\bigr)\subset C^{qr}(\uW)
$$
of bounded cochains, again with a metrizable locally convex topology.
Thus $\delta^{qr}$ restricts to a homomorphism $C_b^{qr}(\uW)\to C_b^{q+1,r}(\uW)$ (but $\overline\partial^{qr}$ is not continuous from $C_b^{qr}(\uW)$ to $C_b^{q,r+1}(\uW))$.
Finally, we let $Z_b^{qr}(\uW)= Z^{qr}(\uW)\cap C_b^{qr} (\uW)$
and endow it with the topology inherited from $C_b^{qr} (\uW)$.

\proclaim{Lemma 2.3}If $\uW$ is a special cover of $M$, $r\geq 1$, and $W_0,\ldots,W_q\in\uW$, then there is a homomorphism
$$
T=T_{W_0\ldots W_q}\colon\Omega_b^{0,r}\bigl(\bigcap^q_{i=0} W_i,E\bigr)
\to\Omega_b^{0,r-1}\bigl(\bigcap^q_{i=0} W_i,E\bigr)\tag2.4
$$
such that $\overline\partial Tf=f$ if $\overline\partial f=0$.
\endproclaim

\demo{Proof}Because of our hypotheses, we can assume that $M=\bC^n$ and $E$ is the trivial bundle $A\times\bC^n\to\bC^n$.
The key is a result of Range and Siu, who construct, when $A=\bC$, homomorphisms $T$ as in (2.4) with the property that $\overline\partial T+T\overline\partial=\id$, see [RS, (3.9)].
In fact, their operator $T$ is a (locally uniform) limit of integral operators $T^\delta$, see [RS, (2.8)---(2.9)].
Now the integrals defining $T^\delta$ make sense if instead of $\bC$--valued forms we substitute in them $A$--valued forms, and by the Banach--Hahn theorem they satisfy the same estimates as the scalar valued operators.
It follows that the $A$--valued integral operators converge to a homomorphism $T$ claimed by the Lemma.
\enddemo

A simple consequence is 

\proclaim{Theorem 2.4}If $M$ is Stein and $E=A\times M\to M$ is a trivial Fr\'echet bundle, then the Dolbeault cohomology groups $H^{0,r}_{\overline\partial} (M,E)$ vanish for $r\geq 1$.
\endproclaim

The only surprising thing about this theorem is that it cannot be generalized much.
Vogt in [V] shows that $H^{0,1}_{\overline\partial} (M,E)$ need not vanish even for $M=\bC$ but with a trivial bundle $E$ with general complete locally 
convex fibers, and Leiterer explained to me that nontrivial Fr\'echet
bundles over Stein manifolds may have nonvanishing 
$H^{0,1}_{\overline\partial}$.
Such Fr\'echet bundles can be constructed from locally trivial fiber bundles with Stein base and fibers, whose total space is not Stein (``counterexamples to Serre's conjecture,'' see [Sk]).

\demo{Proof}First assume that $M$ is an open subset of
some $\Bbb C^k$ and exhaust it by holomorphically convex compact sets 
$K_1,K_2,\ldots$, with $K_i\subset\text{ int }K_{i+1}$.
It follows from Lemma 2.3 that if $f$ is a smooth $\overline\partial$--closed $E$--valued $(0,r)$ form in a neighborhood of $K_i$, then on a possibly smaller neighborhood $f$ is $\overline\partial$--exact, provided $r\geq 1$.
We shall also need the following approximation result, true for any $r\geq 0$:\ the above $f$ can be approximated by $\overline\partial$--closed forms in $\Omega^{0,r}(M,E)$, the approximation taking place in any seminorm $\|\ \ \|_{K_{i-1},p,j}$, see (2.3).

Indeed, if $r\geq 1$, this follows by writing $f=\overline\partial g$ and approximating $g$ in $\|\ \ \|_{K_{i-1},p,j+1}$ by forms 
$g'\in\Omega^{0,r-1}(M,E)$, so that $f'=\overline\partial g'$ will be 
the required approximation of $f$.
If $r=0$, we are talking about approximating holomorphic functions, 
and the result is a rather simple special case of [Lm2, Th\'eor\`eme 1.1].

Let now $r\geq 1$ and $f\in\Omega^{0,r}(M,E)$ be $\overline\partial$--closed.
Fix seminorms $p_1\leq p_2\leq\ldots$ on $A$ defining its topology.
In a neighborhood of $K_i$ write $f=\overline\partial g_i$.
We can arrange recursively that
$$
\|g_{i+1}-g_i\|_{K_{i-1},p_i,i}<2^{-i},\qquad i=1,2,\ldots\tag2.5
$$
Indeed, we note that in general $\overline\partial (g_{i+1}-g_i)=0$ in a neighborhood of $K_i$, hence there is a closed $h\in\Omega^{0,r-1}(M,E)$ such that
$$
\|g_{i+1}-g_i-h\|_{K_{i-1},p_i,i} < 2^{-i}.
$$
Therefore replacing $g_{i+1}$ by $g_{i+1}-h$ we achieve (2.5); which in turn implies that $g=\lim g_i$ is in $\Omega^{0,r-1}(M,E)$ and satisfies $\overline\partial g=f$. Thus we are done if $M$ is open in $\Bbb C^k$.

The general case then follows, since a connected Stein manifold
$M$ can be embedded in a Stein domain $M'\subset\Bbb C^k$ as a holomorphic
retract, and the holomorphic retraction induces a monomorphism
$H^{0,r}_{\overline\partial}(M,E)\to H^{0,r}_{\overline\partial}(M',E')$,
where $E'$ is the trivial bundle $A\times M'\to M'$.
\enddemo

Now we return to a general complex manifold and a Fr\'echet bundle $E\to M$.
If $\uU,\uV$ are open covers of $M$ for which we defined the spaces 
$C^{qr}(\uU), C^{qr}(\uV)$, and $\uU\to\uV$ is a refinement map, then
restriction induces homomorphisms $C^{qr}(\uV)\to C^{qr}(\uU)$.
The image of $f\in C^{qr}(\uV)$ will be denoted $f|\uU\in C^{qr}(\uU)$.

\proclaim{Lemma 2.5}If $\uV$ is a special cover of $M$, $\uU$ a countable refinement of $\uV$ consisting of Stein open sets, and $q,r\geq 0$, then there is a homomorphism
$$
\var=\var^{qr}\colon Z^{qr}(\uU)\to Z_b^{qr}(\uV)
$$
with the following properties:

(a) if $f\in Z^{0r}(\uU)$ then $\var f|\uU=f$;

(b) if $q\geq 1$ and $f\in Z^{qr}(\uU)$, then there is an $h\in C^{q-1,r}(\uU)$ such that
$$
\overline\partial h=0\qquad\text{ and }\qquad\var f|\uU-f=\delta h;
$$

(c) the restriction of $\var$ to $Z_b^{qr}(\uU)$ is continuous;

(d) if $\uU$ is also special and $f\in Z_b^{qr}(\uU)$, then $h$ in (b) can be chosen in $C_b^{q-1,r}(\uU)$.
\endproclaim

Note that (c) is not automatic, as the topology of $Z_b^{qr}(\uU)$ is inherited from $C_b^{qr}(\uU)$, and is different from the one inherited 
from $Z^{qr}(\uU)\subset C^{qr}(\uU)$ when $r>0$.

\demo{Proof}Define homomorphisms 
$$
R=R^{qr}\colon Z^{qr}(\uU)\to C^{q-1,r}(\uU),\qquad S=S^{qr}\colon Z_b^{qr} (\uV)\to C_b^{q,r-1}(\uV)
$$
as follows.
Let $\{\chi_U\colon U\in\uU\}$ be a smooth partition of unity subordinate to $\uU$; then for $f=(f_{U_0\ldots U_q})\in Z^{qr}(\uU)$, $q\geq 1$,
$$
(Rf)_{U_1\ldots U_q}=\sum_{U\in\uU}\chi_U f_{UU_1\ldots U_q}.\tag2.6
$$
If $g=(g_{V_0\ldots V_q})\in Z_b^{qr}(\uV)$, $r\geq 1$, apply the $\overline\partial$ solution operator of Lemma 2.3 componentwise to construct $Sg=(T_{V_0\ldots V_q} g_{V_0\ldots V_q})$.
Thus $\delta Rf=f$ and $\overline\partial Sg=g$.

We prove the Lemma by induction on $q$.
The base case $q=0$ holds since both $Z^{0r}(\uU)$ and $Z_b^{0r}(\uV)$ are naturally identified with the space of $\overline\partial$--closed elements of $\Omega^{0,r}(M,E)$, the topology coming from $Z^{0r}(\uU)$ being finer than the topology coming from $Z_b^{0r}(\uV)$.
If $q\geq 1$, assume that $\var^{q-1,r}$ has been constructed for every $r$, and let
$$
\var^{qr}=\delta S^{q-1,r+1}\var^{q-1,r+1}\overline\partial R^{qr}.
$$
To see that this is indeed a homomorphism $Z^{qr}(\uU)\to Z_b^{qr}(\uV)$, notice that $\delta\overline\partial R=\overline\partial\delta R=0$ and $\overline\partial\delta S=\delta\overline\partial S=0$, so that $\overline\partial R^{qr}$ maps $Z^{qr}(\uU)$ in $Z^{q-1,r+1}(\uU)$ and $\delta S^{q-1,r+1}$ maps $Z_b^{q-1,r+1}(\uV)$ in $Z_b^{qr}(\uV)$.
If $f\in Z^{qr}(\uU)$, let
$$
h_1=(S^{q-1,r+1}\var^{q-1,r+1}\overline\partial R^{qr} f)|\uU-R^{qr} f\in C^{q-1,r} (\uU).\tag2.7
$$
Then $\delta h_1=\var^{qr} f|\uU-f$, and by the inductive assumption
$$
\overline\partial h_1=(\var^{q-1,r+1}\overline\partial R^{qr} f)|\uU-\overline\partial R^{qr} f
$$
either vanishes (if $q=1$), and we take $h=h_1$ in (b); or else there is 
$h_2\in C^{q-2,r+1}(\uU)$ such that $\overline\partial h_2=0$ and 
$\overline\partial h_1=\delta h_2$.
By Theorem 2.4 there is an $h_3\in C^{q-2,r}(\uU)$ with $\overline\partial h_3=h_2$.
Then $h=h_1-\delta h_3$ satisfies 
$\overline\partial h=\overline\partial h_1-\delta\overline\partial h_3=0$ and $\delta h=\var^{qr} f|\uU-f$, which completes the proof of (b).

Next (c) follows since $\overline\partial R$ is continuous from $Z_b^{qr}(\uU)$ to $Z^{q-1,r+1}(\uU)$, as one checks by applying $\overline\partial$ to (2.6).
Finally, if $\uU$ is also special and $f\in Z_b^{qr}(\uU)$, $q\geq 1$, then by induction in (2.7) $h_1\in C_b^{q-1,r}(\uU)$ and $h_2$ above is in $C_b^{q-2,r+1}(\uU)$.
By Lemma 2.3 one can choose $h_3\in C_b^{q-2,r}(\uU)$, whence $h=h_1-\delta h_3\in C_b^{q-1,r}(\uU)$, which proves (d).
\enddemo

\demo{Proof of Theorem 2.2}First we show that (2.2) and (2.1) are monomorphisms.
Being a special cover is stable under small $C^\infty$ perturbation of the boundaries of the covering sets.
Therefore we can construct another special cover $\uV$ such that the 
closure of each $W\in\uW$ is contained in some $V\in\uV$.
Suppose $f\in Z_b^q(\uW,E)=Z_b^{q0}(\uW)$ represents a cohomology class 
$[f]\in H_b^q(\uW,E)$ that (2.2) sends to $0\in H^q(M,E)$.
We apply Lemma 2.5a,b with $\uU=\uW$, to conclude that the image of $[\var^{q0} f]\in H^q (\uV,E)$ in $H^q(M,E)$ factors through $[f]\in H_b^q (\uW,E)$; hence this image is zero.
On the other hand Theorem 2.4 and Dolbeault isomorphism shows that $\uV$ is a Leray cover for the bundle $E$, so that $H^q(\uV,E)\to H^q(M,E)$ is a bijection.
Hence $[\var^{q0} f]=0$, i.e., $\var^{q0} f=\delta g$ with some $g\in C^{q-1,0} (\uV)$.
Since $g|\uW\in C_b^{q-1,0} (\uW)$ and $\var^{q0} f|\uW=\delta (g|\uW)$, Lemma 2.5a,d imply $f$ is the coboundary of a bounded cochain:\ $[f]=0$.
Thus (2.2) is indeed a monomorphism, and with a Stein refinement
$\uU$ of $\uW$ so is (2.1), since the canonical homomorphism 
$H^q(\uU,E)\to H^q(M,E)$ is a bijection, for the same reason as 
$H^q(\uV,E)\to H^q (M,E)$ was.

Next we take an arbitrary Stein refinement $\uU$ of $\uW$ and apply Lemma 2.5 with $\uV=\uW$.
Suppose $f\in Z^q(\uU,E)=Z^{q0}(\uU)$ is $\delta$--exact.
Again, Lemma 2.5a,b imply that the image of $[\var^{q0}f]\in H_b^q (\uW,E)$ in $H^q(M,E)$ is 0; hence by what has been proved, $[\var^{q0} f]=0$.
Thus $\var^{q0}$ maps exact cocycles in $Z^{q0}(\uU)$ to exact cocycles in $Z_b^{q0} (\uW)$.
It induces a homomorphism $$H^q(\uU,E)\to H_b^q(\uW,E),$$
which, by Lemma 2.5a,b is a right inverse to the restriction homomorphism in (2.1).
Passing to the direct limit we obtain a right inverse $H^q(M,E)\to H_b^q (\uW,E)$ of (2.2).
Since (2.1), (2.2) are monomorphisms, the right inverses are left inverses as well, q.e.d.
\enddemo

\proclaim{Corollary 2.6}If $\uW$ is a special cover of $M$ and $\uU$ is a Stein refinement of $\uW$, then
$$
H^q (M,E)\approx H_b^q (\uW,E)\approx H^q (\uU,E).
$$
\endproclaim

\proclaim{Corollary 2.7}If $M$ is compact and $E\to M$ is a holomorphic Banach bundle, then the cohomology groups $H^q(M,E)$ are complete, seminormable topological vector spaces.
\endproclaim

\demo{Proof}Let $\uW$ be a finite special cover of $M$.
Then $H_b^q(\uW,E)$ is a quotient of a Banach space, hence
complete and seminormable; it is also isomorphic to $H^q(M,E)$.
\enddemo

\head 3.\ $\sigma$--compact sets in Fr\'echet spaces\endhead

Recall that a topological space is $\sigma$--compact if it is the countable union of compact subspaces.
A homomorphism of Fr\'echet spaces is compact if the image of some nonempty open set has compact closure.

\proclaim{Lemma 3.1}The range of a compact homomorphism of Fr\'echet spaces is contained in a $\sigma$--compact set.
Conversely, any $\sigma$--compact subset $K$ of a Fr\'echet space $Z$ is contained in the range of a compact homomorphism $Y\to Z$, where $Y$ is a Banach space.
\endproclaim

\demo{Proof}The first part being obvious, we prove the converse only.
Assume first that $K$ is compact; by passing to its hull, we may as well assume $K$ is convex and balanced.
Such a set defines a Banach space $Y=Z_K=\bigcup \{\lambda K\colon\lambda\geq 0\}$, endowed with the norm $\|z\|=\min \{\lambda\geq 0\colon z\in\lambda K\}$.
Then the map $\varphi(z)=z$ from $Y$ to $Z$ is a compact homomorphism and $\varphi(Y)\supset K$.
In general $K=\bigcup K_j$ with $K_j\subset Z$ compact, $j\in\bN$.
For each $j$ there are Banach spaces $(Y_j,\|\ \ \|_j)$ and compact homomorphisms
$\varphi_j\colon Y_j\to Z$ such that $\varphi_j (Y_j)\supset K_j$.
Choose seminorms $p_1\leq p_2\leq\ldots$ on $Z$ defining its topology; we can assume that $p_j(\varphi_j(y))\leq \|y\|_j/2^j$ for $y\in Y_j$.
The $l^1$--sum of the $Y_j$'s
$$
Y=\{ y=(y_j)\colon y_j\in Y_j,\ \|y\|=\sum_j \|y_j\|_j < \infty\},
$$
with the norm $\|\ \ \|$ is a Banach space, $\varphi\colon Y\to Z$ given by
$$
\varphi(y)=\sum_j\varphi_j (y_j),\qquad y=(y_j)\in Y,
$$
is compact, and $\varphi(Y)\supset \bigcup K_j=K$.
\enddemo

\proclaim{Lemma 3.2}Let $C,Z$ be Fr\'echet spaces and $\delta\colon C\to Z$ a homomorphism.
If there is a $\sigma$--compact $L\subset Z$ such that $Z=L+\delta C$, then $\delta C$ is closed and $\dim Z/\delta C < \infty$.
\endproclaim

\demo{Proof}By Lemma 3.1 there are a Banach space $Y$ and a compact homomorphism $\varphi\colon Y\to Z$ with $\varphi Y\supset L$.
Define $\psi,\theta\colon Y\oplus C\to Z$ by
$$
\psi(y,c)=\varphi(y),\qquad\theta (y,c)=\delta(c),
$$
so that $\psi+\theta\colon Y\oplus C\to Z$ is onto.
Since $\psi$ is compact, by Schwartz's theorem
[Sc] $\theta(Y\oplus C)=\delta C$ is indeed closed and of finite codimension.
\enddemo

\head 4.\ Compact perturbations\endhead

Let $M$ be a complex manifold and $E,F$ holomorphic Banach bundles over $M$.
If $U\subset M$ is open, we shall say a homomorphism $k\colon E|U\to F|U$ is fiberwise compact if $k|E_x\in\Hom (E_x,F_x)$ is compact for every $x\in U$.

\proclaim{Proposition 4.1}Suppose $U\subset\subset V\subset M$ 
are open and $E|V$, $F|V$ are trivial. If $k\in\Hom (E|V,F|V)$
is fiberwise compact, then the homomorphism
$$
\alpha\colon\Gamma(V,E)\ni f\mapsto kf|U\in\Gamma(U,F)
$$
is compact.
\endproclaim

This would not be true for Fr\'echet bundles, which is essentially the only reason why Theorem 1.2 is restricted to Banach bundles.

\demo{Proof}Let $(A,\|\ \ \|_A), (B,\|\ \ \|_B)$ be Banach spaces such that $E|V$, $F|V$ are (isomorphic to) $A\times U$ resp.~$B\times U$.
The set
$$
\cF=\{f\in\Gamma (V,E)\colon\sup_{x\in U} \|f(x)\|_A < 1\}
$$
is open in $\Gamma(V,E)$.
Its image under $\alpha$ is a locally equicontinuous family on 
$U$ by Cauchy estimates; also $\{\alpha f(x)\colon f\in \cF\}$ 
has compact closure in $B$ for every $x\in U$.
By the Arzel\'a--Ascoli theorem any sequence in $\alpha(\cF)$
has a subsequence convergent in the topology of $\Gamma(U,F)$, 
and so $\alpha$ is indeed compact.
\enddemo

\proclaim{Proposition 4.2}If $F$ is a compact perturbation of $E$,
then $E$ is also a compact perturbation of $F$.
More precisely, if $\uU$ is an open cover of $M$ and 
$\varphi_U\colon E|U\to F|U$ are fredomorphisms such that 
$\varphi_U-\varphi_V$ are fiberwise compact for $U,V\in\uU$, then there 
are an open cover $\uU'$ of $M$ and fredomorphisms 
$\psi_{U'}\colon F|U'\to E|U'$ for $U'\in\uU'$ such that
$\psi_{U'} -\psi_{V'}$, $\psi_{U'}\varphi_U-id_{E|U'\cap U}$, and 
$\varphi_U\psi_{U'}-id_{F|U'\cap U}$ are fiberwise compact, $U\in\uU$, 
$U', V'\in \uU'$.
\endproclaim

\demo{Proof}Given $x\in U\in\uU$, there is a finite codimensional subspace $E'_x\subset E_x$ on which $\varphi_U$ is injective.
In fact, $E'_x$ can be extended to a holomorphic subbundle $E'$ of $E$, over a neighborhood $N(x)\subset U$ of $x$, and $\varphi_U|E'$ is still injective.
It follows that $F'=\varphi_U (E')\subset F|N(x)$ is a holomorphic subbundle of finite corank, and $\varphi_U$ restricts to an isomorphism $\varphi'\colon E'\to F'$.
Upon shrinking $N(x)$ we can assume $F'$ has a complementary bundle $F''\subset F|N(x)$, and in fact, at the price of refining $\uU$ we can assume that the bundles $E',F',F''$ are defined over all of $U=N(x)$.
If $p'\colon F|U=F'\oplus F''\to F'$ denotes the projection, then
$$
\psi_U=\varphi^{\prime -1} p'\in\Hom (F|U, E|U)
$$
is a fredomorphism.
Furthermore $\psi_U\varphi_U-\text{id}_{E|U}=k_U$ and 
$\varphi_U\psi_U-\text{id}_{F|U}=l_U$ are fiberwise compact (in fact, 
of finite rank).
With $\psi_U$, $U\in\uU=\uU'$ thus constructed
$$
\psi_V-\psi_U=\psi_U (\varphi_U-\varphi_V)\psi_V-k_U\psi_V+\psi_U l_V
$$
are also fiberwise compact, whence the claim follows.
\enddemo

From now on we assume $M$ is compact.
It will be convenient to use indexed covers $\uU=\{ U_i\in I\}$,
write $U_{ij\ldots}=U_i\cap U_j\cap\ldots$, and denote the components 
$f_{U_i U_j\ldots}$ of a cochain by $f_{ij\ldots}$.

Let $\uU=\{U_i\colon i\in I\}$ and $\uV=\{V_i\colon i\in I\}$ be finite open covers, $U_i\subset\subset V_i$.
Assume that $E|V_i$ and $F|V_i$ are trivial, and that there are fredomorphisms $\varphi_i\colon E|V_i\to F|V_i$ such that $\varphi_i-\varphi_j$ are fiberwise compact.
Define homomorphisms of $\Phi^q\colon C^q(\uV,E)\to C^q (\uU,F)$ between spaces of holomorphic cochains by
$$
(\Phi^q f)_{i_0\ldots i_q}=\varphi_{i_0} f_{i_0\ldots i_q}|U_{i_0\ldots i_q},\qquad
f\in C^q (\uV,E).\tag4.1
$$

\proclaim{Proposition 4.3}$\Phi^q \delta-\delta\Phi^{q-1}\colon C^{q-1} (\uV,E)\to C^q(\uU,F)$ is compact for $q\geq 1$.
\endproclaim

\demo{Proof}If $f=(f_{i_1\ldots i_q})\in C^{q-1} (\uV,E)$ then
$$
\align
& (\Phi^q\delta f)_{i_0\ldots i_q}=\varphi_{i_0}\sum^q_{j=0} (-1)^j f_{i_0\ldots\hat i_j\ldots i_q}|U_{i_0\ldots i_q} \qquad\text{ and}\\
& (\delta\Phi^{q-1} f)_{i_0\ldots i_q}=\varphi_{i_1} f_{i_1\ldots i_q} +\varphi_{i_0}\sum^q_{j=1} (-1)^j f_{i_0\ldots\hat i_j\ldots i_q}|U_{i_0\ldots i_q}.
\endalign
$$
Hence the claim follows from Proposition 4.1, as
$$
(\Phi^q\delta f-\delta\Phi^{q-1} f)_{i_0\ldots i_q}=(\varphi_{i_0}-\varphi_{i_1}) f_{i_1\ldots i_q}|U_{i_0\ldots i_q}.
$$
\enddemo

\demo{Proof of Theorem 1.2}Construct three finite Stein covers
$\uU=\{U_i\colon i\in I\}$, $\uV=\{V_i\colon i\in I\}$, and $\uW=\{W_i\colon i\in I\}$ of $M$, $U_i\subset\subset V_i\subset\subset W_i$.
Make sure $\uW$ is special for both $E$ and $F$ (see Definition 2.1) and
so fine that there are fredomorphisms $\varphi_i\colon E|W_i\to F|W_i$ and 
$\psi_i\colon F|W_i\to E|W_i$ such that
$$
\varphi_i-\varphi_j,\quad\psi_i-\psi_j,\quad
\varphi_i\psi_j-\id_{F|W_{ij}},\text{ and }\psi_i\varphi_j-\id_{E|W_{ij}}
$$
are all fiberwise compact, see Proposition 4.2.
Define homomorphisms $\Phi^q$ as in (4.1) and $\Psi^q\colon C^q (\uW,F)\to C^q (\uV,E)$ by
$$
(\Psi^q f)_{i_0\ldots i_q}=\psi_{i_0} f_{i_0\ldots i_q}|V_{i_0\ldots i_q},\qquad f\in C^q (\uW,F).
$$
It follows from Proposition 4.1 that
$$
\Phi^q\Psi^q f=f|\uU+\kappa f,\qquad\text{ where }\quad\kappa\colon C^q (\uW,F)\to C^q(\uU,F)\tag4.2
$$
is compact.
By the case $r=0$ of Lemma 2.5 there are homomorphisms
$$
\var=\var^q\colon Z^q(\uU,F)\to Z^q(\uW,F)
$$
such that $\var f|\uU$ is cohomologous to $f$ for every $f\in Z^q (\uU,F)$.
In view of (4.2), then
$$
\Phi^q \Psi^q\var\equiv \id_{Z^q(\uU,F)}+\kappa\var\qquad \mod\delta C^{q-1} (\uU,F).\tag4.3
$$
By Corollary 2.6 $H^{q+1}(\uV,E)\approx H^{q+1} (M,E)$ is Hausdorff, that is,
$\delta C^q(\uV,E)\subset Z^{q+1}(\uV,E)$ is closed.
Michael's selection theorem [M, Theorem 3.2$''$] implies there is a continuous map
$$
\mu\colon\delta C^q (\uV,E)\to C^q (\uV,E)
$$
with $\delta\mu=\id_{\delta C^q(\uV,E)}$.
(Typically, Michael's theorem is formulated for Banach, rather than Fr\'echet spaces.
However, it has always been understood that the proof applies for Fr\'echet spaces as well, see e.g.~the last paragraph on page 364 in [M].)
Consider the continuous map
$$
\beta=\Psi^q\var-\mu\delta\Psi^q\var\colon Z^q (\uU,F)\to C^q (\uV,E).\tag4.4
$$
In fact, $\beta$ maps into $Z^q(\uV,E)$, since $\delta\beta=\delta \Psi^q\var-\delta\mu\delta \Psi^q\var=0$.

Now we come to the main point:\ for every $\sigma$--compact $K\subset Z^q(\uV,E)$ there is a $\sigma$--compact $L\subset Z^q(\uU,F)$ such that
$$
\beta^{-1}(K+\delta C^{q-1} (\uV,E))\subset L+\delta C^{q-1} (\uU,F).\tag4.5
$$
Indeed, by (4.3), (4.4), and since $\delta|Z^q (\uW,F)=0$
$$
\align
\Phi^q\beta&=\Phi^q\Psi^q\var-\Phi^q\mu\delta\Psi^q\var\\
&\equiv \id_{Z^q(\uU,F)}+\kappa\var+\Phi^q\mu (\Psi^{q+1}\delta-\delta\Psi^q)\var\qquad\mod\delta C^{q-1}(\uU,F),
\endalign
$$
cf. (4.2). Here $\kappa$ is a compact operator, and so is
$\Psi^{q+1}\delta-\delta\Psi^q$, see Proposition 4.3.
Hence by Lemma 3.1 for any $S\subset Z^q(\uU,F)$ we have
$$
S\subset\Phi^q \beta S+H+\delta C^{q-1} (\uU,F),\tag4.6
$$
with some $\sigma$--compact $H\subset C^q (\uU,F)$.
If now $S$ stands for the left hand side of (4.5), then
$$
\align
\Phi^q\beta S&=\Phi^q (K+\delta C^{q-1} (\uV,E))\\
&\subset\Phi^q K+(\Phi^q\delta-\delta\Phi^{q-1}) C^{q-1} (\uV,E)+\delta C^{q-1}(\uU,F).
\endalign
$$
The first term in the last line is $\sigma$--compact, and by Lemma 3.1
and Proposition 
4.3  the second is contained in a $\sigma$--compact set.
Combining this with (4.6) we obtain a $\sigma$--compact $L\subset C^q(\uU,F)$ such that
$$
S\subset L+\delta C^{q-1} (\uU,F).
$$
Since $L$ can be replaced by $L\cap Z^q (\uU,F)$, (4.5) is verified.

By hypothesis, $H^q(M,E)\approx H^q(\uV,E)$ is finite dimensional.
We take for $K\subset Z^q (\uV,E)$ a finite dimensional complementary subspace to $\delta C^{q-1} (\uV,E)$.
In this case (4.5) gives
$$
Z^q (\uU,F)=L+\delta C^{q-1} (\uU,F)
$$
with some $\sigma$--compact $L$.
Hence Corollary 2.6 and Lemma 3.2 imply
$$
H^q (M,F)\approx H^q (\uU,F)=Z^q (\uU,F)/\delta C^{q-1} (\uU,F)
$$
is indeed Hausdorff and finite dimensional.
\enddemo

Theorem 1.2 also holds on compact manifolds with strongly pseudoconvex
boundaries, by essentially the same proof, and, as said,
it should have a sheaf version as well. But in the setting of this paper
there does not seem to be much room to improve on it.
Without the assumption that $H^{q+1}(M,E)$ is Hausdorff it would not hold.
The reasoning to follow was inspired by an idea that I learned from Leiterer, who attributed it to V\^ aj\^ aitu.
Suppose $M$ is K\"ahler and $H^1(M,\cO)\neq 0$, so that there is 
a sequence of nontrivial line bundles $\Lambda_k\to M$ 
converging to the trivial line bundle.
If $E$ is (a suitable Hilbertian completion of) $\bigoplus\limits_1^\infty\Lambda_k$, and $F$ the trivial Hilbert bundle $l^2\times M\to M$, then $E,F$ are compact perturbations of one another, $\dim H^0 (M,E)<\infty$ but $\dim H^0(M,F)=\infty$.
In more detail, let $\uU$ be a Stein cover of $M$ and $g=(g_{UV})\in Z^1 (\uU,\cO)$ not exact.
The multiplicative cocycles $(e^{g_{UV}/k})\in Z^1 (\uU,\cO^*)$ are not exact either for large $k$; by rescaling $g$ we can assume they are not exact for any $k\in\bN$.
Thus the line bundles $\Lambda_k$ they determine 
are holomorphically nontrivial.
Since topologically they nevertheless are trivial, $H^0 (M,\Lambda_k)=0$.
The infinite diagonal matrices $h_{UV}$ with diagonal entries $e^{g_{UV}/k}, k\in\bN$, define a cocycle $h=(h_{UV})\in Z^1(\uU,\text{GL}(l^2))$ and so a Hilbert space bundle $E\to M$.
Each $\Lambda_k$ canonically embeds in $E$, $\bigoplus\Lambda_k$ is dense in $E$, and there are holomorphic projections $E\to\Lambda_k$.
It follows that $H^0 (M,E)=0$.
However, the trivial bundle $F=l^2\times M\to M$, a compact perturbation of $E$, has $\dim H^0 (M,F)=\infty$.

By Theorem 1.2 we conclude that $H^1(M,E)$ is not Hausdorff.
This example reveals one more notable thing:\ while all groups $H^q(M,F)$ are Hausdorff, a compact perturbation $E$ has a non Hausdorff cohomology group.---Non Hausdorff cohomology groups of holomorphic Hilbert bundles over compact manifolds were first constructed by Erat in [E].
His construction is similar, being based on a finite rank vector bundle over $\bP_1$ that can be deformed in a non--isomorphic bundle.

\enddocument
\bye

\bye